\documentclass[12pt,a4paper]{amsart}

\usepackage{amssymb,amsmath,amsfonts,amsthm}
\usepackage{graphicx}
\graphicspath{ {./images/} }

\usepackage{mathrsfs}  
\usepackage{psfrag}
\usepackage{mathtools}
\usepackage{color}
\usepackage{todonotes}
\usepackage{enumitem}

\usepackage{hyperref}
\hypersetup{colorlinks=false, citecolor=green, linkcolor=red, hypertexnames=false}

\usepackage{chngcntr}
\counterwithin{equation}{section}

\theoremstyle{plain}
\newtheorem{main}{Theorem}

\newtheorem{maincor}[main]{Corollary}
\newtheorem{theorem}{Theorem}[section]

\theoremstyle{remark}

\newcommand{\Sing}{\operatorname{Sing}}

           \def\ea{\end{array}}
          \def\ec{\end{center}}
     \def\ed{\end{description}}
        \def\ee{\end{equation}}
       \def\eea{\end{eqnarray}}
     \def\eeaa{\end{eqnarray*}}
 \def\et{\end{thebibliography}}

\def\bM{{\bf{M}}}

\def\Orb{{\rm Orb}}

\def\Sing{{\rm Sing}}

\def\cG{{\mathcal G}}

\def\cD{{\mathcal D}}

\def\cU{{\mathcal U}}

\def\cR{{\mathcal R}}

\def\cF{{\mathcal F}}
\def\cM{{\mathcal M}}
\def\cN{{\mathcal N}}
\def\cP{{\mathcal P}}

\def\cR{{\mathcal R}}

\def\cS{{\mathcal S}}

\def\vep{\varepsilon}

\def\RR{{\mathbb R}}

\def\Exp{\operatorname{Exp}}

\title[Equilibrium states for sectional-hyperbolic flows]{Uniqueness of equilibrium states for Lorenz attractors in any dimension}
\date{\today}
\author{Maria Jose Pacifico, Fan Yang and Jiagang Yang}

\address{Instituto de Matem\'atica, Universidade Federal do Rio de Janeiro, C. P. 68.530, CEP 21.945-970,  Rio de Janeiro, RJ, Brazil.}
 \email{pacifico@im.ufrj.br }

\address{Department of Mathematics, Michigan State University, East Lansing, Michigan, USA.}
\email{yangfa31@msu.edu}

\address{Departamento de Geometria, Instituto de Matem\'atica e Estat\'\i stica, Universidade Federal Fluminense, Niter\'oi, Brazil}
\email{yangjg\@@impa.br}

\thanks{This research has been supported [in part] by CAPES - Finance Code $001$ and CNPq-grants. MJP was partially supported by  by Grant Cientista do Nosso Estado (FAPERJ)}

\begin{document}

\maketitle

\begin{abstract}
	In this note we consider the thermodynamic formalism for Lorenz attractors of flows in any dimension. Under a mild condition on the H\"older continuous potential function $\phi$, we prove that for an open and dense subset of $C^1$ vector fields, every Lorenz attractor supports a unique equilibrium state. In particular, we obtain the uniqueness for the measure of maximal entropy.
\end{abstract}

\section{Introduction and statement}
Ever since its discovery in 1963 by Lorenz~\cite{Lo63}, the {\em Lorenz attractor} has been playing a central role in the research of singular flows, i.e., flows generated by smooth vector fields with {\em singularities}. The investigation of their statistical and topological properties has been successful in dimension three, thanks to the geometric Lorenz model~\cite{GW79}; however, results in higher dimensions are surprisingly lacking. 

This note is devoted to the study of the statistical properties of Lorenz attractors in any dimension.
{In particular,  we are concerned with {\em equilibrium states} for H\"older continuous potential functions $\phi:\bM\to \RR$ on Lorenz attractors. }

In what follows, $\bM$ is a compact Riemannian manifold without boundary, and $\mathscr{X}^1(\bM)$ denotes the space of $C^1$ vector fields on $\bM$ endowed with $C^1$ topology.  Given $X\in \mathscr{X}^1(\bM)$, a  singularity is a point $q\in\bM$ where $X(q)=0$. The set of singularities of $X$ is denoted by $\Sing(X)$.

An $X$-invariant set $\Lambda$ is called a {\em Lorenz attractor} if it satisfies the following conditions ($\varphi_t$ denotes the flow generated by $X$):
\begin{enumerate}
	\item $\varphi_t\mid\Lambda$ is {\em sectional-hyperbolic}, meaning that there exists an $D\varphi_t$-invariant dominated splitting $E^s\oplus E^{cu}$ such that $D\varphi_t\mid E^s$ is uniformly contracting, and $D\varphi_t\mid E^{cu}$ is {\em volume expanding}: there exists $C>0,\lambda>0$ such that 
	$$
	|\det D\phi_t(x)\mid_{V_x}| \ge C e^{\lambda t} 
	$$
	for every $x\in\bM$, every subspace $V_x\subset E_x^{cu}$ with $\dim V_x\ge 2$, and every $t>0$;
	\item $\varphi_t\mid\Lambda$ is {\em transitive}: there exists a point $x\in\Lambda$ whose orbit is dense in $\Lambda$;
	\item  there exists a neighborhood $U\supset \Lambda$ such that for every $x\in\Lambda$ one has $\omega(x)\subset \Lambda$. 
\end{enumerate} 

Recall that given such a function $\phi$, the variational principle for flows states that
$$
P(\phi) = \sup_{\mu\in\cM(X)} \left(h_\mu(\varphi_1) + \int\phi\,d\mu\right),
$$
where $P(\phi)$ is the topological pressure of $\phi$, $\cM(X)$ is the space of $X$-invariant probability measures, and $h_\mu(\varphi_1)$ is the measure-theoretic entropy of the time-one map $\varphi_1$. An equilibrium state is an invariant  measure $\mu\in\cM(X)$ that achieves the supremum. As a special case, the topological pressure of the constant function $\phi\equiv 0$ is the topological entropy of the flow, and any equilibrium state in this case is called a {\em measure of maximal entropy}.

It has been proven in~\cite{PYY} that for every sectional-hyperbolic set $\Lambda$ of a flow $\varphi_t$, $\varphi_t\mid\Lambda$ is {\em entropy expansive}. Due to the celebrated work of Bowen~\cite{B72}, the metric entropy, as a function on $\cM(X)$, is upper semi-continuous. Consequently, every continuous function $\phi:\bM\to\RR$ admits an equilibrium state $\mu_\phi$. Our main theorem below states that under certain mild conditions on $X$ and $\phi$, this equilibrium state is unique.

\begin{main}\label{m.A}
	There exists an open and dense subset $\cR\subset\mathscr{X}^1(\bM)$, such that for every $X\in\cR$ and every Lorenz attractor $\Lambda$ of $X$, let $\phi:\bM\to\RR$ be a H\"older continuous function such that 
	$$
	\phi(\sigma) < P(\phi), \forall \sigma\in\Sing(X\mid\Lambda).
	$$ 
	Then there exists a unique equilibrium state $\mu$ supported in $\Lambda$.
\end{main}

As a corollary, we obtain
\begin{maincor}
	There exists an open and dense subset $\cR\subset\mathscr{X}^1(\bM)$, such that for every $X\in\cR$ and every Lorenz attractor $\Lambda$ of $X$, there exists a unique measure of maximal entropy for $X\mid \Lambda$.
\end{maincor}

The corollary easily follows from Theorem~\ref{m.A} applied to the constant function $\phi\equiv 0$ (note that every Lorenz attractor has positive topological entropy, by~\cite[Theorem C]{PYY}).

\section{A revised Climenhaga-Thompson criterion for the uniqueness of equilibrium states}
The main tool for the proof of Theorem~\refeq{m.A} is the following theorem from~\cite{PYY21}, which improves the Climenhaga-Thompson criterion~\cite{CT16} by weakening the specification assumption. 

We define
\begin{equation*}\label{e.L}
	L_X = \max_{t\in[0,1]} L_{\varphi_t}\ge 1
\end{equation*}
where $L_{\varphi_t}$ is the Lipschitz constant of $\varphi_t$. Writing
$$
\Gamma_\vep(x)=\{y\in \bM: d(\varphi_t(x),\varphi_t(y))\le \vep \mbox{ for all }t\in\RR\},
$$
for the bi-infinite Bowen ball at $x$, a vector field $X$ is said to be {\em almost expansive} at scale $\vep>0$, if the set 
\begin{equation*}\label{e.Exp}
	\Exp_\vep(X) :=  \left\{x\in \bM: \Gamma_\vep(x)\subset \varphi_{[-s,s]}(x) \mbox{ for some }s>0\right\}
\end{equation*}
has full probability: for any $\mu\in\cM(\bM)$ one has $\mu(\Exp_\vep(X))=1$. Clearly, almost expansivity is weaker than expansivity; however, it also implies entropy expansivity at the same scale~\cite{LVY}.

Below we will identify a pair $(x,t)\in\bM\times\RR^+$ with the orbit segment $\{\varphi_s(x):0\le s\le t\}$. Given a collection of orbit segments $\cD\subset \bM\times\RR^+$, a {\em decomposition}
$(\cP,\cG,\cS)$ of $\cD$ consists of three collections $\cP,\cG,\cS\subset \bM\times \RR^+$ and three functions $p,g,s:\cD\to \RR^+$ such that for every $(x,t)\in \cD$, the values $p=p(x,t), g=g(x,t)$ and $s=s(x,t)$ satisfy $t=p+g+s$, and
\begin{equation}\label{e.decomp}
	(x,p)\in\cP,\hspace{0.5cm} (\varphi_p(x),g)\in\cG,\hspace{0.5cm} (\varphi_{p+g}(x),s)\in\cS.
\end{equation}
Now we are ready to state the main result of~\cite{PYY21}.

\begin{theorem}\label{t.criterion}
		Let $(\varphi_t)_{t\in\RR}$ be a Lipschitz continuous flow on a compact metric space $\bM$, and $\phi:\bM\to\RR$ a continuous function. Suppose that there exist $\vep>0,\delta>0$ with $\vep\ge1000L_X\delta$ such that $X$ is almost expansive at scale $\vep$, and there exists $\cD\subset\bM\times\RR^+$ which admits a decomposition $(\cP,\cG,\cS)$ with the following properties:
	\begin{enumerate}[label={(\Roman*)}]
		\item $\cG$ has tail (W)-specification at scale $\delta$;
		\item  $\phi$ has the Bowen property at scale $\vep$ on $\cG$;
		\item $P(\cD^c\cup [\cP]\cup[\cS], \phi,\delta,\vep)<P(\phi)$.
	\end{enumerate}
	Then there exists a unique equilibrium state for the potential $\phi$.
\end{theorem}

We invite the interested readers to~\cite{CT16,PYY21} for the precise definition of the involved terms. 

\section{Structure of the proof}

From now on, $\Lambda$ is a Lorenz attractor unless otherwise specified. We will assume that all the singularities of $X\mid\Lambda$ are hyperbolic; this is clearly an open and dense condition among $C^1$ vector fields. 

We will prove Theorem~\refeq{m.A} by verifying the assumptions of Theorem~\ref{t.criterion}. 
We first remark that Theorem~\ref{t.criterion} can be applied, with only minor modification if necessary, to $X\mid \Lambda$. 
Also note that the almost expansivity at some given scale $\delta_0>0$ has been proven in~\cite[Theorem F]{PYY21}. This implies almost expansivity at any smaller scale.
We will first construct two orbit segment collections $\cG_\delta$ and $\cG_0$ on which the specification and the Bowen property hold, respectively (Theorem~\ref{t.s},~\ref{t.B}).  $\cG_\delta$ and $\cG_0$ must be taken large enough, such that one could find $\cG\subset \cG_\delta\cap\cG_0$ with ``large topological pressure'' (Theorem~\ref{t.decomposition}).

\subsection{The weak specification property}

The next theorem establishes the specification property on a collection of orbit segments. 

\begin{theorem}\label{t.s}
	For every $\delta>0$, there exists an orbit collection $\cG_\delta\subset \bM\times\RR^+$ such that $\cG_\delta$ has the tail (W)-specification property at scale $\delta$.
\end{theorem}

We remark that the choice of the orbit segment $\cG_\delta$ depends on the scale $\delta$. It is worth noting that in the original Climenhaga-Thompson criterion~\cite[Theorem 2.9]{CT16}, the specification assumption is made on $\cG^M$ for every $M\ge 0$ (at the same scale $\delta$), where $\cG^M$ is defined as (see~\eqref{e.decomp})
$$
\cG^M = \{(x,t)\in\cD:p\vee s\le M \}.
$$
Thinking of $\cG$ as the `good core' of $\cD$, it is clear that
$$
\cG^M\subset \cG^{M'} \mbox{ whenever } M<M', \mbox{ and } \bigcup_{M\ge0}\cG^M = \cD.
$$
In the subsequent applications of~\cite{CT16}, this assumption is usually 
replaced by the assumption that $\cG$ has the tail (W)-specification property at {\em all scales}. Unfortunately in our case, we do not know whether $\cap_{\delta>0}\cG_\delta$ carries large pressure. To deal with this issue, in~\cite{PYY21} we relaxed the specification condition of~\cite{CT16} and assume that it only holds on the `good core' $\cG$ (Theorem~\ref{t.criterion} (I)).

The proof of Theorem~\ref{t.s} calls for the following topological description for Lorenz attractor, which strengthens the main result of~\cite{CY}.

An invariant set $\Lambda$ is called {\em Lyapunov stable} if for every neighborhood $U$ of $\Lambda$, there exists a neighborhood $V$ of $\Lambda$, such that for every $x\in V$ we have $\varphi_t(x)\in U, \forall t>0$. 

\begin{theorem}\label{t.t}
	There exists a $C^1$ residual subset $\cR\subset \mathscr{X}^1(\bM)$, such that for every $X\in\cR$ and every sectional-hyperbolic, Lyapunov stable chain recurrent class $\Lambda$ of $X$, there exists a neighborhood $\cU$  of $X$ and a neighborhood $U$ of $\Lambda$ that satisfy the following properties:
	\begin{enumerate}
		\item every $Y\in\cU$ has a unique chain recurrent class $C_Y\subset U$ which is a sectional-hyperbolic attractor;
		\item $C_Y$ is the homoclinic class of some hyperbolic periodic orbit $p_Y\in C_Y$; in particular, $W^s(\Orb(p_Y))$ is dense in $C_Y$;
		\item for every point $x\in C_Y\setminus\Sing(Y)$, $\cF_Y^s(x)\pitchfork W_Y^u(p_Y)\ne\emptyset$.
	\end{enumerate}
\end{theorem}

Here $\cF_Y^s(x)$ is the leaf of the stable foliation $\cF_Y^s$ for the flow $Y$ that contains $x$. It is well known (see, for instance, \cite[Lemma 3.10]{PYY}) that $\Lambda\cap\left(W^{ss}(\sigma)\setminus \{\sigma\}\right) =\emptyset$ for every sectional-hyperbolic set\footnote{Recall that all the singularities in a sectional-hyperbolic set are {\em Lorenz-like}, meaning that $E^s_\sigma =  E_\sigma^{ss}\oplus_< E_\sigma^c$ with $\dim E_\sigma^c = 1$; furthermore, the bundle $E^s$ from the dominated splitting $E^s\oplus E^{cu}$ on $\Lambda$ matches the bundle $E^{ss}_\sigma$ at the singularity $\sigma$, and $\cF^s_Y(\sigma)$ coincides with the strong stable manifold $W^{ss}(\sigma)$ tangent to $E^{ss}_\sigma$ at $\sigma$.} $\Lambda$ and every $\sigma\in\Lambda\cap\Sing(X)$, so (3) is optimal. 


\subsection{The Bowen property}
We prove the following theorem:
\begin{theorem}\label{t.B}
	There exists $\vep_0>0$ and a collection of orbit segments $\cG_0\subset \bM\times \RR^+$, such that $\phi$ has the Bowen property on $\cG_0$ at scale $\vep_0$ (and consequently, at every smaller scale).
\end{theorem}

The proof heavily relies on the scaled linear Poincar\'e flow of Liao~\cite{Liao96}. We write 
\begin{equation*}\label{e.tubular}
	B_{\beta} (x,t) = \bigcup_{s\in[0,t]} N_{X,\varphi_s(x)}(\beta |X(\varphi_s(x))|)
\end{equation*}
where $N_{X,y}(r)$ is the image of the $r$-neighborhood of $0_y$ in the normal plane  $\cN_{X,y} = \langle X(y)\rangle ^\perp\subset T_y\bM$ under the exponential map $\exp_y$. One should think of $B_{\beta} (x,t)$ as a tubular neighborhood of the orbit segment $(x,t)$ whose size is proportional to the flow speed $|X(\varphi_s(t))|$. 
We prove that there exists $\beta>0$ such that every $(x,t)\in\cG_0$ satisfies
$$
\varphi_s(y)\in B_{\beta} (x,t), \forall s\in[0,t] \mbox{ whenever }y\in B_{\vep_0,t}(x)
$$
where $B_{\vep_0,t}(x)$ is the $(\vep_0,t)$-Bowen ball of $x$. Then the Bowen property follows whenever the terminal point $\varphi_t(x)$ is a $\lambda_0$-quasi hyperbolic point for some prescribed $\lambda_0\in(0,1)$ (See~\cite{Liao} and~\cite[Definition 11]{PYY}; also note that vectors in the stable bundle $E^s$ are uniformly contracted by the scaled linear Poincar\'e flow $\psi^*_t$).

\subsection{The decomposition and the pressure gap}
Let $\vep_0>0$ and $\cG_0\subset\bM\times\RR^+$ be given by Theorem~\ref{t.B}, and $\delta_0$ be the scale of the almost expansivity given by~\cite[Theorem F]{PYY}. We take $\vep <\min\{\vep_0,\delta_0\}$ small enough\footnote{The smallness of $\vep$ is needed to control the pressure function $P(\cdot, \phi,\delta,\vep)$ with respect to its second scale.}, and fix $\delta<\vep/(1000L_X)$. By theorem~\refeq{t.s}, we obtain the orbit collection $\cG_\delta$ which has the specification property. The next theorem shows the existence of an orbit collection $\cD$ with a decomposition.

\begin{theorem}\label{t.decomposition}
	There exists an orbit segments collection $\cD\subset \bM\times\RR^+$ which admits a decomposition $(\cP,\cG,\cS)$ with the following properties:
	\begin{enumerate}
		\item $\cG\subset \cG_0\cap\cG_\delta$;
		\item orbit segments in $\cD^c$, $\cP$ and $\cS$ must spend a significant proportional of their time in a prescribed, small neighborhood $U_\Sing$ of $\Sing(X)\cap\Lambda$.
	\end{enumerate}
\end{theorem}

Property (1) shows that $\cG$ has the Bowen property at scale $\vep$ (Theorem~\refeq{t.B}) and has the tail (W)-specification property at scale $\delta$ (Theorem~\ref{t.s}). 

To show that the pressure of $\cD^c\cup\cP\cup\cS$  is strictly smaller than $P(\phi)$, we used a modified version of the variational principle which shows that the topological pressure of $\cD^c\cup\cP\cup\cS$ is bounded from above by the metric pressure $P_\mu(\phi):= h_\mu(X) + \int\phi\,d\mu $, where $\mu$ is a limit point of some convex combination of empirical measures supported on orbit segments in $\cD^c\cup\cP\cup\cS$. It follows from  Property (2) that $\mu(U_\Sing)>a>0$ for some constant $a$. Using a semi-continuity argument which is made possible by~\cite[Theorem A]{PYY}, we show that there exists $b>0$ such that 
$$
P_\mu(\phi)<P(\phi)-b.
$$
This verifies Assumption (III) of Theorem~\ref{t.criterion}, and finishes the proof of Theorem~\ref{m.A}.

\end{document}